\documentclass[11pt]{amsart}

\usepackage[T1]{fontenc}
\usepackage[utf8]{inputenc}
\usepackage{amsmath,amssymb,amsthm}
\usepackage{enumitem}
\usepackage[a4paper,margin=3cm]{geometry}
\usepackage{placeins}
\usepackage{tikz}
\usepackage[hidelinks,pdfencoding=auto]{hyperref}

\newtheorem{theorem}{Theorem}[section]
\newtheorem{proposition}[theorem]{Proposition}

\theoremstyle{definition}

\title[A note on lower bounds for numerical series]
{A note on lower bounds for numerical series}

\author{R. \'Alvarez-Nodarse}
\address{Departamento de Análisis Matemático, Universidad de Sevilla}
\address{c/Tarfia s/n, 41012 Sevilla, Spain}
\email{ran@us.es}

\author{K. Castillo}
\address{CMUC, Department of Mathematics, University of Coimbra,
3000-143 Coimbra, Portugal}
\email{kenier@mat.uc.pt}

\date{29 July 2026}

\hypersetup{
  pdftitle={A note on lower bounds for numerical series},
  pdfauthor={R. Alvarez-Nodarse and K. Castillo},
  pdfsubject={A mathematical correction and sharp reciprocal-moment
  inequalities},
  pdfkeywords={moment sequences, reciprocal moments,
  positive linear functionals, sharp lower bounds,
  Hurwitz zeta function}
}

\begin{document}

\begin{abstract}
After passage to the representing measure, the three principal
theorems of Esstafa and Sfaxi, \emph{J.\ Math.\ Anal.\ Appl.}
\textbf{556} (2026), 130199, reduce respectively to Jensen's inequality,
one quadratic identity, and a one-point measure.  We solve the
underlying reciprocal-moment problem sharply, determine all
extremisers, and obtain a strictly stronger Hurwitz-zeta bound.  We also
identify the gap between formal inversion and evaluation at \(1\), and
give explicit counterexamples to further analytic claims.
\end{abstract}

\maketitle

\section{Introduction}

Esstafa and Sfaxi~\cite{ES26} study lower bounds for a convergent moment
series associated with a normalised positive linear functional.  The
Hamburger representation turns their assumptions into the single
identity
\[
 S:=\sum_{n\ge0}a_n
 =\mathbb E\!\left[\frac1{1-X}\right],\quad
 -1<X<1,\quad\text{almost surely}.
\]
It contains the entire substance of their three principal theorems.
The universal bound follows from the stronger moment-dependent Jensen
bound; every value in \((1/2,\infty)\) is realised by a one-point
measure; and the estimate involving the first two moments follows from
one quadratic identity.  None requires the formal-inverse or
Jacobi-functional machinery developed in~\cite{ES26}.

The same reduction raises the natural question not addressed there:
determine the best lower bound for \(\mathbb E[Y^{-1}]\) when the range,
mean, and variance of \(Y\) are prescribed.  We obtain the exact bound,
identify all extremal laws, and derive a strict improvement of the
published Hurwitz-zeta estimate.

The algebraic operations used in~\cite{ES26} are purely formal; they
supply neither analytic convergence nor boundary evaluation.  The
resulting gap, together with further false analytic assertions, is
recorded at the end.

\section{The elementary content of the published theorems}

Let \((a_n)_{n\ge0}\) be a real sequence with \(a_0=1\), and let
\(\mathbf u\) be the linear functional on \(\mathbb R[x]\) determined by
\(\mathbf u(x^n)=a_n\).  The standing assumptions in~\cite{ES26} are
\begin{equation}\label{eq:hypotheses}
 \sum_{n\ge0}a_n\ \text{converges},\quad
 \mathbf u(p^2)\ge0,\quad p\in\mathbb R[x].
\end{equation}
By the Hamburger moment theorem~\cite{ST43}, there is a probability
measure \(\mu\) on \(\mathbb R\) such that
\[
 a_n=\int_{\mathbb R}x^n\,\mathrm d\mu(x),\quad n\ge0.
\]

\begin{proposition}\label{prop:reduction}
Under~\eqref{eq:hypotheses}, the measure \(\mu\) is carried by
\((-1,1)\), and
\begin{equation}\label{eq:sum-integral}
 \sum_{n\ge0}a_n
 =
 \int_{(-1,1)}\frac{\mathrm d\mu(x)}{1-x}
 <\infty.
\end{equation}
\end{proposition}

\begin{proof}
Convergence gives \(a_{2n}\to0\), while
\[
 a_{2n}
 =\int_{\mathbb R}x^{2n}\,\mathrm d\mu(x)
 \ge\mu\bigl(\{x:|x|\ge1\}\bigr).
\]
Thus \(\mu\) is carried by \((-1,1)\).  Pairing consecutive terms and
using monotone convergence,
\begin{align*}
 \sum_{n\ge0}a_n
 &=\lim_{N\to\infty}
   \int_{(-1,1)}
   (1+x)\sum_{j=0}^{N}x^{2j}\,\mathrm d\mu(x)\\[7pt]
 &=\int_{(-1,1)}\frac{\mathrm d\mu(x)}{1-x}.
\end{align*}
\end{proof}

\subsubsection*{Theorem~1 of~\cite{ES26} follows immediately from
Jensen's inequality}

Let \(X\) have distribution \(\mu\), and write
\[
 S:=\sum_{n\ge0}a_n
 =\mathbb E\!\left[\frac1{1-X}\right].
\]
Since \(-1<X<1\) almost surely,
\[
 S-\frac12
 =\mathbb E\!\left[\frac{1+X}{2(1-X)}\right]>0.
\]
Although \(1/2\) is the optimal universal constant, convexity gives the
stronger moment-dependent estimate
\[
 S\ge\frac1{1-a_1}>\frac12,
\]
with equality in the first inequality precisely when
\(a_2-a_1^2=0\).

\subsubsection*{Theorem~3 of~\cite{ES26} is realised by a one-point
measure}

Given \(T>1/2\), put \(c:=1-T^{-1}\).  Then \(-1<c<1\);
\(\delta_c\) is a probability measure, hence gives a normalised
positive functional, and its moments \(c^n\) are absolutely summable
with sum \((1-c)^{-1}=T\).  Thus every
\(T>1/2\) occurs.  Conversely, the preceding identity gives \(S>1/2\)
for every admissible sequence; the limiting choice \(c=-1\) is
inadmissible because \(\sum_{n\ge0}(-1)^n\) diverges.  Hence the range is
\[
 \left(\frac12,\infty\right).
\]
Thus \(1/2\) is the infimum, not a minimum, and no Jacobi functional is
required.

\subsubsection*{Theorem~2 of~\cite{ES26} follows from one identity}

Put
\[
 m:=\mathbb E[X]=a_1,\quad
 v:=\operatorname{Var}(X)=a_2-a_1^2,\quad
 Y:=1-X,
\]
and suppose that \(v>0\).  Since \(0<Y<2\), if \(r:=S^{-1}\), then
\begin{align*}
 1-m-r
 &=\mathbb E\!\left[\frac{(Y-r)^2}{Y}\right]\\[7pt]
 &>\frac12\,\mathbb E[(Y-r)^2]\\[7pt]
 &\ge\frac v2.
\end{align*}
The first inequality is strict because \(v>0\).
It follows simultaneously that
\[
 1-m-\frac v2>0
\]
and that
\begin{equation}\label{eq:published-bound}
 S>\frac1{1-m-\frac12v}.
\end{equation}

\section{The sharp reciprocal-moment bound}

For non-zero variance, the published estimate is not sharp.  The
resulting two-moment extremal problem is classical in form; its exact
solution for the present data follows from a quadratic minorant with
two contacts.

\begin{theorem}\label{thm:sharp-reciprocal}
Let \(Y\) be a random variable such that
\[
 0<Y\le L,\quad\text{almost surely},
\]
where \(L>0\).  If \(\mathbb E[Y]=L\), then \(Y=L\) almost surely and
\(\mathbb E[Y^{-1}]=L^{-1}\).  Otherwise, put
\[
 s:=\mathbb E[Y],\quad v:=\operatorname{Var}(Y).
\]
Then
\begin{equation}\label{eq:feasible}
 0\le v<s(L-s)
\end{equation}
and
\begin{equation}\label{eq:sharp-general}
 \mathbb E\!\left[\frac1Y\right]
 \ge
 \frac{L^2-Ls-v}{L(Ls-s^2-v)}.
\end{equation}
The expectation on the left may be \(+\infty\).  If
\[
 h:=s-\frac{v}{L-s},
\]
then equality holds precisely when the law of \(Y\) is carried by
\(\{h,L\}\).  Thus the bound is attained for every pair \((s,v)\)
satisfying~\eqref{eq:feasible}.  If \(Y<L\) almost surely, equality is
possible only when \(v=0\); for \(v>0\) the same bound is optimal but is
not attained.
\end{theorem}

\begin{proof}
Since \(s<L\), the event \(\{Y<L\}\) has positive probability.  Hence
\[
 Ls-s^2-v=\mathbb E[Y(L-Y)]>0.
\]
Set
\[
 h:=s-\frac{v}{L-s}>0
\]
and
\[
 q_h(y):=
 \frac1L+\frac{2(L-y)}{Lh}
 +\frac{y^2-Ly}{Lh^2}.
\]
The identity
\begin{equation}\label{eq:minorant-identity}
 \frac1y-q_h(y)
 =
 \frac{(y-h)^2(L-y)}{Lh^2y},
 \quad 0<y\le L,
\end{equation}
shows that \(q_h(y)\le y^{-1}\).  Since
\(Ls-s^2-v=h(L-s)\),
\begin{align*}
 \mathbb E[q_h(Y)]
 &=\frac1L+\frac{2(L-s)}{Lh}
   +\frac{s^2+v-Ls}{Lh^2}\\[7pt]
 &=\frac{L^2-Ls-v}{L(Ls-s^2-v)}.
\end{align*}
This proves~\eqref{eq:sharp-general}.  The right-hand side of
\eqref{eq:minorant-identity} vanishes only at \(h\) and \(L\), which
proves the equality assertion.  Conversely, the two-point law
\[
 \mathbb P(Y=h)=\frac{L-s}{L-h},\quad
 \mathbb P(Y=L)=\frac{s-h}{L-h}
\]
has mean \(s\), variance \((s-h)(L-s)=v\), and attains equality.

It remains only to verify optimality under the strict condition
\(Y<L\).  When \(v>0\), \eqref{eq:feasible} gives
\(s+v/s<L\).  Choose \(b\in(s+v/s,L)\), and set
\[
 h_b:=s-\frac{v}{b-s}>0,
\]
and let \(Y_b\) take the values \(h_b\) and \(b\) with probabilities
\[
 \frac{b-s}{b-h_b}
 \quad\text{and}\quad
 \frac{s-h_b}{b-h_b},
\]
respectively.  Then
\[
 \mathbb E[Y_b]=s,\quad
 \operatorname{Var}(Y_b)=v,
\]
while
\[
 \mathbb E[Y_b^{-1}]
 =\frac{b+h_b-s}{bh_b}
 \longrightarrow
 \frac{L^2-Ls-v}{L(Ls-s^2-v)}
\]
as \(b\uparrow L\).
\end{proof}

Figure~\ref{fig:minorant} displays the two contacts which determine
both the equality cases and the limiting extremisers.

\begin{figure}[htbp]
\centering
\begin{tikzpicture}[
  x=3.1cm,
  y=1.25cm,
  line cap=round,
  line join=round
]
  \tikzset{
    axis/.style={draw=black!70,line width=0.68pt,-stealth},
    reciprocal/.style={draw=black!82,line width=1.0pt},
    minorant/.style={
      draw=black!48,
      line width=0.9pt,
      dash pattern=on 4pt off 2.4pt
    },
    guide/.style={
      draw=black!30,
      line width=0.5pt,
      dash pattern=on 2pt off 2pt
    },
    label/.style={font=\scriptsize,text=black!72,inner sep=1.5pt}
  }

  \path[use as bounding box] (0.12,-0.22) rectangle (2.22,3.58);

  \begin{scope}
    \clip (0.12,-0.22) rectangle (2.22,3.58);
    \clip
      plot[domain=0.10:2,samples=150]
        (\x,{1/\x})
      --
      plot[domain=2:0.10,samples=150]
        (\x,{0.5+(2-\x)/0.8+(\x*\x-2*\x)/(2*0.8*0.8)})
      -- cycle;
    \shade[left color=white,right color=black!18]
      (0.18,-0.22) rectangle (0.62,3.80);
    \fill[black!18]
      (0.62,-0.22) rectangle (2.01,3.80);
  \end{scope}

  \begin{scope}
    \clip (0.20,0) rectangle (2.04,3.58);
    \draw[reciprocal]
      plot[domain=0.10:2,samples=150]
        (\x,{1/\x});
    \draw[minorant]
      plot[domain=0.10:2,samples=150]
        (\x,{0.5+(2-\x)/0.8+(\x*\x-2*\x)/(2*0.8*0.8)});
  \end{scope}

  \draw[axis] (0.20,0) -- (2.15,0)
    node[right,label] {$y$};
  \draw[axis] (0.20,0) -- (0.20,3.45);

  \draw[guide] (0.8,0) -- (0.8,1.25);
  \draw[guide] (2,0) -- (2,0.5);

  \fill[black!78] (0.8,1.25) circle (1.45pt);
  \fill[black!78] (2,0.5) circle (1.45pt);

  \node[label,below] at (0.8,-0.02) {$h$};
  \node[label,below] at (2,-0.02) {$L$};

  \draw[reciprocal] (1.20,2.82) -- (1.40,2.82);
  \node[label,anchor=west] at (1.45,2.82) {$y^{-1}$};

  \draw[minorant] (1.20,2.50) -- (1.40,2.50);
  \node[label,anchor=west] at (1.45,2.50) {$q_h(y)$};
\end{tikzpicture}
\caption{The quadratic minorant in
\eqref{eq:minorant-identity}, shown for \(L=2\) and \(h=4/5\).
It is tangent to \(y^{-1}\) at \(h\) and meets it again at \(L\).
The shaded gap is the non-negative remainder in
\eqref{eq:minorant-identity}; the two contacts identify all extremal
distributions.  The diffuse left edge marks the artificial truncation
of the plot as \(y\downarrow0\).}
\label{fig:minorant}
\end{figure}

For the series in~\eqref{eq:hypotheses}, take \(Y=1-X\) and \(L=2\).
Theorem~\ref{thm:sharp-reciprocal} gives
\begin{equation}\label{eq:sharp-main}
 S\ge
 \frac{2+2m-v}{2(1-m^2-v)}.
\end{equation}
For \(v>0\), the inequality in~\eqref{eq:sharp-main} is strict because
\(Y<2\) almost surely.  Moreover, by~\eqref{eq:feasible},
\[
 1-m-\frac v2
 >
 1-m-\frac{1-m^2}{2}
 =\frac{(1-m)^2}{2}>0,
\]
so every denominator below is positive.  The right-hand side
of~\eqref{eq:sharp-main} is strictly greater than the published bound,
since
\begin{equation}\label{eq:strict-improvement}
 \frac{2+2m-v}{2(1-m^2-v)}
 -\frac1{1-m-\frac12v}
 =
 \frac{v^2}
 {2(2-2m-v)(1-m^2-v)}>0.
\end{equation}
The approximating two-point laws in the proof, with \(L=2\), correspond
under \(X=1-Y\) to finite-support measures in \((-1,1)\); their moment
series converge absolutely.  Thus~\eqref{eq:sharp-main} is sharp within
the original admissible class.
Consequently, for \(v>0\) the complete hierarchy is
\[
 S
 >
 \frac{2+2m-v}{2(1-m^2-v)}
 >
 \frac1{1-m-\frac12v}
 >
 \frac1{1-m}
 >
 \frac12.
\]

\subsection{A sharper Hurwitz-zeta bound}

Let \(\alpha>0\), \(s>1\), and set
\[
 A:=\left(\frac{\alpha}{\alpha+1}\right)^s,\quad
 B:=\left(\frac{\alpha}{\alpha+2}\right)^s.
\]
Since \(\alpha(\alpha+2)<(\alpha+1)^2\),
\[
 0<A^2<B<A<1.
\]
The probability density
\[
 \frac{\alpha^s}{\Gamma(s)}
 x^{\alpha-1}(-\log x)^{s-1},
 \quad 0<x<1,
\]
has moments
\[
 \mathbb E[X^n]
 =\left(\frac{\alpha}{n+\alpha}\right)^s,
\]
and hence
\[
 \mathbb E\!\left[\frac1{1-X}\right]
 =\alpha^s\zeta(s,\alpha).
\]
Applying Theorem~\ref{thm:sharp-reciprocal} to \(Y=1-X\) with \(L=1\)
gives
\begin{equation}\label{eq:hurwitz-sharp}
 \alpha^s\zeta(s,\alpha)
 >
 \frac{A^2+A-B}{A-B}.
\end{equation}
Writing \(v=B-A^2\), the quantity denoted by \(X(s,\alpha)\)
in~\cite{ES26} is \(A+v/2\).  Consequently, bound~(13) there is
\[
 \alpha^s\zeta(s,\alpha)>\frac1{1-A-v/2}.
\]
Here \(v>0\), \(A(1-A)-v=A-B>0\), and
\(2-2A-v>(1-A)(2-A)>0\).  The improvement
in~\eqref{eq:hurwitz-sharp} is therefore strict because
\[
 \frac{A^2+A-B}{A-B}
 -\frac1{1-A-v/2}
 =
 \frac{v(A+v)}
 {(2-2A-v)\bigl(A(1-A)-v\bigr)}>0.
\]

\FloatBarrier

\section{Further comments}

\subsubsection*{The formal-inverse argument}

Esstafa and Sfaxi~\cite{ES26} use the algebraic dual \(\mathcal P'\);
no topology on \(\mathcal P\) is specified.  This is legitimate for
their algebraic operations.  Writing
\((\mathbf v)_n:=\mathbf v(x^n)\), the correspondence
\[
 M_{\mathbf v}(z):=\sum_{n\ge0}(\mathbf v)_n z^n
\]
identifies the Cauchy product in \(\mathcal P'\) with multiplication in
the ring \(\mathbb R[[z]]\) of formal power series.  In particular,
\[
 M_{\mathbf u^{-1}}(z)M_{\mathbf u}(z)=1,\quad
 \text{in }\mathbb R[[z]].
\]
This is a coefficientwise identity.

This formal identity has no automatic boundary value at \(z=1\).
In the proof of Theorem~2 in~\cite[p.~12]{ES26}, the functional
\[
 w:=\frac{-x^2\mathbf u^{-1}}{a_2-a_1^2}
\]
is shown to be normalised and positive, but the convergence of
\(\sum_{n\ge0}(w)_n\), required by Theorem~1, is not established.
They subsequently use the evaluation
\begin{equation}\label{eq:formal-evaluation}
 \sum_{n\ge0}(\mathbf u^{-1})_n
 =
 \left(\sum_{n\ge0}a_n\right)^{-1}.
\end{equation}
This is not a formal consequence: it requires convergence of the
inverse series and an Abelian passage to \(z=1\).  No such argument is
given in~\cite{ES26}.  Placing \(\mathcal P\) or \(\mathcal P'\) in a
topological-vector-space setting would not by itself supply the missing
convergence.  Consequently, the published argument does not prove
Theorem~2.  The theorem itself remains true for the elementary reason
given in Section~2.

\subsubsection*{The Laurent expansion}

On pp.~6--7, \cite{ES26} writes
\[
 F(z)=\sum_{n\ge0}\frac{a_n}{z^{n+1}}
     =\int\frac{\mathrm d\mu(x)}{z-x},
 \quad z\in\mathbb C\setminus(-1,1).
\]
The integral is analytic off the topological support of \(\mu\).
The displayed Laurent series converges absolutely and represents the
integral for \(|z|>1\).  This does not justify the asserted identity
throughout \(\mathbb C\setminus(-1,1)\): at points with \(|z|<1\) the
series may diverge, while boundary points require separate analysis.
For the probability measure with density \(1\) on \((-1,0)\),
\[
 a_n=\frac{(-1)^n}{n+1},\quad
 \sum_{n\ge0}a_n=\log2.
\]
At \(z=i/2\),
\[
 \left|\frac{a_n}{z^{n+1}}\right|
 =\frac{2^{n+1}}{n+1}\not\longrightarrow0,
\]
although the integral is well defined, whereas at \(z=-1\) the
integral itself diverges.  The displayed equality therefore fails on
its stated domain.

The half-plane mapping asserted in the same passage is also reversed.
If \(\operatorname{Im}z>0\), then
\[
 \operatorname{Im}F(z)
 =
 -\operatorname{Im}z
 \int_{(-1,1)}\frac{\mathrm d\mu(x)}{|z-x|^2}<0,
\]
whereas~\cite{ES26} gives the opposite sign and concludes that \(F\)
maps the upper half-plane into itself.  In fact \(F\) maps it into the
lower half-plane; with the standard convention, \(-F\) is the
Nevanlinna function.

\subsubsection*{The Riemann-zeta application}

On p.~18, \cite{ES26} defines
\[
 X(s):=
 2^{-s}+\frac{1}{2\cdot3^s}
 \left(1-\left(\frac34\right)^s\right)
\]
and asserts that \(X(s)>1/2\) for every \(s>1\).  At \(s=2\),
\[
 X(2)
 =\frac14+\frac1{18}\left(1-\frac9{16}\right)
 =\frac{79}{288}<\frac12.
\]
The resulting claim \(1/\zeta(s)<1/2\) is also false, since
\[
 \frac1{\zeta(2)}=\frac6{\pi^2}>\frac12.
\]
The conclusion that the functional \(v_s\), with moments
\[
 (v_s)_n=\frac{\mu_{\mathrm M}(n+1)}{(n+1)^s},
\]
is not positive is correct, but not for the reason given.  Here
\(\mu_{\mathrm M}\) is the Möbius function, and
\[
 H_1(v_s)
 =
 \begin{vmatrix}
  1&-2^{-s}\\[7pt]
  -2^{-s}&-3^{-s}
 \end{vmatrix}
 =-3^{-s}-4^{-s}<0.
\]

\subsubsection*{Funding declaration}

This work was supported by the Portuguese Foundation for Science and
Technology (FCT) under the project UID/00324/2025, Centre for
Mathematics of the University of Coimbra.  RAN was partially supported
by PID2024-155593NB-C21 (FEDER(EU)/Ministerio de Ciencia, Innovación y
Universidades--Agencia Estatal de Investigación), IMUS-Maria de Maeztu grant CEX2024-001517-M - 
Apoyo a Unidades de Excelencia María de Maeztu funded by MICIU/AEI/ DOI: 10.13039/501100011033,
and FQM-415 (Junta de Andalucía). 
KC  was supported by FCT project 2022.00143. CEECIND/CP1714/CT0002
(\href{https://doi.org/10.54499/2022.00143.CEECIND/CP1714/CT0002}
{DOI: 2022.00143.CEECIND}).

\end{document}